\documentclass[12pt]{article}
\usepackage{amssymb}
\usepackage{amsmath}
\usepackage[latin1]{inputenc}
\usepackage[T1]{fontenc}
\usepackage[dvips]{graphicx}

\def \G{\mathbb{G}}

\def \N{\mathbb{N}}

\def \Q{\mathbb{Q}}

\def \Z{\mathbb{Z}}

\newcommand{\ord}{\operatorname{ord}}

\newenvironment{demo}[1]{{\bf Proof: }}{\hfill\mbox{$\Box$} \vskip 0.7cm}

\newtheorem{theorem}{Theorem}[section]
\newtheorem{proposition}[theorem]{Proposition}

\newtheorem{corollary}[theorem]{Corollary}
\newtheorem{lemma}[theorem]{Lemma}
\newtheorem{example}[theorem]{Example}

\title{Computing the torsion points of a variety defined by lacunary polynomials}
\author{Louis Leroux}
\date{}
\begin{document}

\maketitle
\thispagestyle{empty}
\begin{abstract}
\noindent We present an algorithm for computing the set of torsion points satisfying a given system of multivariate polynomial equations. Its comple\-xity is quasilinear in the logarithm of the degree of the input equations and exponential in their number of non zero terms and variables.
\vskip 0.5cm
\noindent{\bf 2010 MS Classification:} Primary $11$Y$16$, Secondary $12$Y$05$, $68$W$30$.
\end{abstract}

\let\thefootnote\relax\footnotetext{The author was partially supported by the CNRS PICS "Properties of heights of arithmetic varieties",  $(2009-2011)$.}

\section{Introduction and statement of results}

Let $F_1,\ldots,F_k\in\Z[X_1,\ldots,X_n]$ be a family of multivariate polynomials. We consider the problem of computing the solutions of the system of equations
\begin{equation}
\label{equation1}
F_1(X_1,\ldots,X_n)=\cdots=F_k(X_1,\ldots,X_n)=0
\end{equation}
in  roots of unity. This problem arises naturally when solving trigonome\-tric equations and also when solving some geometric problems, see for instan\-ce  \cite{CJ76}, \cite{PR98}.

For $n=1$, this is equivalent to the problem of finding the common cyclotomic factors of $F_1,\ldots,F_k$. Based on a result of J.H. Conway and A.J. Jones \cite{CJ76}, M. Filaseta, A. Granville and A. Schinzel have recently described an algorithm in \cite{FGS08}, that solves this problem with a complexity quasilinear in the logarithm of the degree of the input polynomials.

For $n\geq 2$, the system $(\ref{equation1})$ can have an infinite number of solutions in roots of unity but a structure theorem of M. Laurent \cite{Lau84} implies that this set can in principle be described in finite terms (see below for more details). In this text, we first simplify the algorithm of \cite{FGS08} and we precise its complexity (see Algorithm $1$ in Section $2$), then we extend it to the general multivariate case (see Algorithm $2$ in Section $3$). The complexity of this algorithm is again quasilinear in the logarithm of the degree of the input polynomials, although exponential in their number of variables and their number of non zero terms. Hence our algorithm can be regard as an effective version of Laurent's Theorem which is particulary well suited for lacunary polynomials, that is polynomials having a small number of non zero terms but potentially large degree.

\bigskip

In the following, we fix $n\in\N^*$. We begin by recalling some basic facts on Diophantine geometry of subvarieties of the algebraic torus $\G_m^n:=(\overline{\Q}^*)^n$. All definitions and results stated here can also be found in \cite{S96} or in \cite{Z09} and we refer to these texts for the proofs and for more details. The set $\G_m^n$ is a group for the usual multiplication:
$$
(x_1,\ldots,x_n)\cdot(y_1,\ldots,y_n)=(x_1y_1,\ldots,x_ny_n).
$$
An {\it algebraic subgroup} $H\subset\G_m^n$ is a subgroup of $\G_m^n$ which is an algebraic variety, i.e. it is Zariski closed. A point of finite order in the group $\G_m^n$ will be called a {\it torsion point}. The set of torsion points in $\G_m^n$ is exactly $\mu_{\infty}^n$, where $\mu_{\infty}$ denotes the set of roots of unity in $\overline{\Q}$. A {\it torsion coset} is a coset of the form
$$
\underline{\eta}H:=\left\{\underline{\eta}\cdot\underline{h}\mid\underline{h}\in H\right\},
$$
where $H$ is an algebraic subgroup and $\underline{\eta}\in\mu_{\infty}^n$.
For $V$ a subvariety of $\G_m^n$, we set
$$
V_{\rm tors}:=V\cap \mu_{\infty}^n,
$$
for the set of torsion points lying in $V$.
A celebrated theorem of Laurent~\cite{Lau84} asserts that the Zariski closure $\overline{V_{\rm tors}}$ of $V_{\rm tors}$ is the union of a finite number of torsion cosets:
\begin{equation}
\label{eqLaurent}
\overline{V_{\rm tors}}=B_1\cup\cdots\cup B_t.
\end{equation}
This result was previously conjectured by S. Lang \cite{Lan83}. It was proved by Y. Ihara, J.P. Serre and J. Tate when $V$ is a curve \cite{Lan83} and by Laurent in the general case.

\bigskip

For $\underline{x}=(x_1,\ldots,x_n)$ $\in\G_m^n$ and $\underline{\lambda}=(\lambda_1,\ldots,\lambda_n)\in\Z^n$, we set:
$$
\underline{x}^{\underline{\lambda}}:=x_1^{\lambda_1}\cdots x_n^{\lambda_n}\in\G_m.
$$
If $\Lambda$ is a subgroup of $\Z^n$ of dimension $k$, then the set 
$$
H_{\Lambda}:=\left\{\underline{x}\in\G_m^n\mid \, \underline{x}^{\underline{\lambda}}=1,\,\,\forall\, \underline{\lambda}\in\Lambda\right\}
$$
is an algebraic subgroup of $\G_m^n$ of dimension $n-k$. Furthermore the map $\Lambda\mapsto H_{\Lambda}$ is a bijection between the set of subgroups of $\Z^n$ and the set of algebraic subgroups of $\G_m^n$ \cite[Lemma 2]{S96}.

Let $M_{k,n}(\Z)$ be the set of matrices with $k$ rows and $n$ columns with integer coefficients. When $L=\left(\lambda_{i,j}\right)_{{1\leq i\leq k}\atop{1\leq j\leq n}}\in M_{k,n}(\Z)$, where $\underline{x}=(x_1,\ldots,x_n)\in\G_m^n$ we define: 
$$
\underline{x}^L:=\left(\underline{x}^{\underline{\lambda}_1},\ldots,\underline{x}^{\underline{\lambda}_k} \right)\in \G_m^k,
$$
where $\underline{\lambda_1},\ldots,\underline{\lambda_k}$ are the rows of the matrix $L$. Now let $\Lambda$ be a subgroup of $\Z^n$ of dimension $k$, let $\underline{\lambda_1},\ldots,\underline{\lambda_k}$ be a $\Z-$basis of $\Lambda$ and let $L$ be the matrix of size $k\times n$ whose rows consist of the vectors $\underline{\lambda_1},\ldots,\underline{\lambda_k}$. Let also $\underline{\omega}\in\mu_{\infty}^k$, then the set:
\begin{equation}
\label{eqoutput}
B(L,\underline{\omega}):=\left\{ \underline{x}\in\G_m^n\mid \underline{x}^{L}=\underline{\omega}\right\}
\end{equation}
is a torsion coset and, in fact, all torsion cosets can be described in this way.
Our algorithm takes as input a finite number of polynomials and outputs a finite number of torsion cosets whose union is $\overline{V_{\rm tors}}$ as in (\ref{eqLaurent}). Furthermore each torsion coset of dimension $n-r$ will be represented by a pair $(L,\underline{\omega})\in M_{k,n}(\Z)\times \mu_{\infty}^k$ as in (\ref{eqoutput}). We remark that it is easy to describe all the torsion points lying in such a coset since
$$
B(L,\underline{\omega})_{\rm tors}:=\left\{ \underline{\zeta}\in\mu_{\infty}^n\mid \underline{\zeta}^{L}=\underline{\omega}\right\}.
$$

\bigskip

Let us now describe the input of the algorithm. Let 
$$
F(X_1,\ldots,X_n)=\displaystyle\sum_{i=1}^Na_i\underline{X}^{\underline{\alpha_i}}
$$
be a polynomial in $n$ variables with integer coefficients. The {\it lacunary enco\-ding} of $F$ is the list
$$
\left[(a_1,\underline{\alpha_1}),\ldots,(a_N,\underline{\alpha_N})\right]
$$
of its non zero coefficients with corresponding exponents and its {\it height} is defined by:
$$
{\rm h}(F):=\max_{1\leq i\leq N}\{\log |a_i|\}.
$$
If $F$ is of total degree $d$ and of height $h$, the number of bits to encode this representation is
$$
O\big(N(h+n\log(d))\big).
$$
When $d\in\N$, we set $M(d)$ for the complexity of multiplying two integers smaller than $d$. We can now state our main result.
\begin{theorem}
\label{variete}
There exists a deterministic algorithm which has the following propriety. Let $F_1,\ldots,F_k\in\Z[X_1,\ldots,X_n]$ be polynomials given by their lacunary encoding and let $V$ be the variety defined by those polynomials. The algorithm computes a family of torsion cosets $B_1,\ldots,B_t$ such that
$$
\bigcup_{1\leq i\leq t} B_i=\overline{V_{\rm tors}},
$$
where the torsion cosets are represented as in $(\ref{eqoutput})$.
Furthermore, this algorithm performs at most
$$
O\left(N^{nkN}(M(d)\log\log(d)+h)\right)
$$
bit operations, where $N$ denotes the maximum number of non zero terms in each $F_i$, $d$ denotes an upper bound for their total degree and $h=\displaystyle{\max_{1\leq i\leq k}}\left({\rm h}(F_i)\right)$.
\end{theorem}
{\bf Remark}: We have
$$
M(d)=O\left(\log(d)\log\log(d)\log\log\log(d)\right)
$$
because of the algorithm of A. Schönnhage and V. Strassen (\cite{GG03}) and so the complexity of the algorithm underlying Theorem \ref{variete} is
$$
O_{\varepsilon}\left(N^{nkN}(\log(d)^{1+\varepsilon}+h)\right) \hskip 1cm \forall\,\varepsilon>0.
$$

\vskip 0.5cm

Some related problems have been studied in the univariate case. Let $F\in\Z[X]$, given by its lacunary encoding, and let $m$ be a non-negative integer. The problem of testing whether $F$ vanishes at a root of unity of order $m$ is called the {\it cyclotomic test} (CT). In \cite{CTV09}, Q. Cheng, S. P. Tarasov and M. N. Vyalyi have shown that this problem can be solved in time polynomial in the size of the input, that is CT$\in$P.

The {\it generalized cyclotomic test} (GCT) consists in determining if $F$ va\-nishes in some root of unity at all, that is if there exists some $m\in\N$ such that $(F,m)$ is a positive instance of CT. The fact that CT$\in$P implies that GCT$\in$NP. On the other hand, Plaisted has previously shown that this pro\-blem is NP-hard \cite[Thm 5.1]{P84} and so GCT is NP-complete.

Filaseta and Schinzel proposed a subexponential algorithm for solving GCT \cite{FS04}. In collaboration with Granville, they extend this algorithm to compute all of the cyclotomic factors of $F$ \cite{FGS08}. In this text, we simplify this last algorithm and we clarify the dependance of its complexity in the number $N$ of terms, which turns out to be exponential. We also improve the practical implementation and time execution of this algorithm. The exponential behavior of the complexity with respect to $N$ seems unavoidable and indeed, we think that the size of the output is exponential in the size of the input in the worst case. In Section $2$, we exhibit a family of examples which support this conjecture.

For $n\geq 2$, the only previoius constructive methods are described in \cite{Ru93}, \cite{BS02}, \cite{AS08}, \cite{Ro07}. Their complexity is not explicited in these references but it is certainly at least of type $d^{2^n}h$, where $n$ denotes their number of variables, $d$ denotes a bound for the degree of the input polynomials, and $h$ denotes a bound for their height, because of the systematic application of iterated resultants. Thus these methods are better suited for dense polynomials of small degree and having a small number of variables. On the contrary, our result is particulary efficient for polynomials with few terms but potentially very high degree ($10^{100}$ for example), see for instances the examples in Subsection~\ref{exemple}.

\vskip 0.5cm

We deduce from the analysis of the algorithm underlying Theorem \ref{variete} the following upper bound for the number of torsion cosets in a subvariety of $\G_m^n$:
\begin{corollary}
\label{cor}
Let $F_1,\ldots,F_k\in\Z[X_1,\ldots,X_n]$ be polynomials, each having at most $N$ non zero terms and let $V$ be the subvariety of $\G_m^n$ defined by those polynomials.
Then $\overline{V_{\rm tors}}$ is the union of at most 
$$
(N!)^k\exp\left({3(n+1)\sqrt{kN\log kN}}\right)
$$
torsion cosets.
\end{corollary}

Historically, the first effective upper bound for the number of maximal torsion cosets has been obtained by E. Bombieri and U. Zannier in \cite{BZ95}. Shortly afterwards, W. Schmidt found an upper bound depending only on the number of variables and on the degree of the input polynomials \cite{S96}. 
The main interest of the bound in Corollary \ref{cor} is that it is independant of the degree of the polynomials defining $V$. However, we mention that a similar bound might be alternatively obtained with the methods in \cite{S96}.

It is worth noting that for the number of connected (or irreducible) torsion cosets, the dependance on the degree is unavoidable. Currently, the best bound is due to F. Beukers and C.J. Smyth \cite{BS02} when $n=2$ and to F. Amoroso and E. Viada in the general case \cite{AV09}. The result of Amoroso and Viada says that the number of connected torsion cosets in a variety defined by polynomials of degree at most $d$ is bounded by
$$
d^n\left(200n^5\log(n^2d)\right)^{n^2(n-1)^2}.
$$
Other results in this direction were obtained by G. Rémond in \cite{R02}, Amo\-ro\-so and S. David in \cite{AD06} and David and P. Philippon in \cite{DP07}.

\vskip 0.5cm

The outline of the paper is as follows. In Section $2$, we present our simplification of the algorithm of Filaseta, Granville and Schinzel. In Section $3$, we generalize it for a multivariate polynomial and we will finally explain in Section \ref{Sec4} how we can adapt it for general varieties. 

\vskip 0.5cm

\noindent{\bf Acknowledgements:} I would like to thank Francesco Amoroso and Martín Sombra for their precious help and many suggestions. I also express my gratitude to Denis Simon who suggested me the example in Section~\ref{exemple}.

\section{Cyclotomic factors of univariate polynomials}

\label{Sec2}

Let $F(X)\in\Z[X]$ be a polynomial of degree $d$. Since $F$ has integer coefficients, computing the set of $\zeta\in\mu_{\infty}$ such that $F(\zeta)=0$ is equivalent to computing the set of integers $m$ for which
$$
\Phi_m(X)|F(X),
$$
where $\Phi_m$ denotes the $m$th cyclotomic polynomial. We want to compute this set in time polynomial in $\log(d)$. However, the total number of cyclotomic factors of $F$ in the worst case is not polynomial in $\log(d)$, as shown by the following example: let $x$ be an integer greater than $17$ and let
$$
F(X)=X^{d}-1,
$$
where $d$ is the product $p_1\cdots p_r$ of all prime numbers smaller than $x$. An effective version of Tchebychev's Theorem \cite{RS62} gives the inequalities:
$$
\pi(x)\geq \frac{x}{\log x}\hskip 0.5cm \text{and}\hskip 0.5cm \sum_{p\leq x}\log p\leq 1.02\,x\hskip 1cm \text{for }x\geq 17,
$$
where $\pi(x):=\#\{p \text{ prime }|\,\,p\leq x\}$. From these, we deduce:
$$
r\geq \dfrac{x}{\log x}\hskip 1cm\text{and}\hskip 0.5cm d\leq \exp(1.02\,x).
$$
Furthermore,
$$
F(X)=\prod_{m|d}\Phi_m(X) \hskip 0.5cm \text{and}\hskip 0.5cm \#\{m\in\N:m|d\}=2^r\geq 2^{\frac{x}{\log x}},
$$
so the number of cyclotomic factors of $F$ is not polynomial in the logarithm of its degree.

Even in the case when the number of factors is small, separating them is as difficult as factorizing integers. To illustrate this problem, let us consider 
$$
F(X):=X^{pq}-1,
$$
where $p$ and $q$ are distinct prime numbers. We have
$$
F(X)=\Phi_1(X)\Phi_p(X)\Phi_q(X)\Phi_{pq}(X)
$$
so we want the algorithm to output the set $\{1,p,q,pq\}$. Doing this is equi\-valent to factor $pq$ but this problem is known to be hard and, at the moment, no algorithm can factor $pq$ in polynomial time in $\log (pq)$.

To avoid these problems, we will represent the output of our algorithm differently. It will be given as a set $S_F$ of pairs of integers $(m,e)$ such that
$$
V(F)_{\rm tors}=\bigcup_{(m,e)\in S_F}V(\Phi_m(X^e)),
$$
where $V(F):=\{x\in\overline{\Q}:F(x)=0\}$.
Given this representation, we are reduced to integer factorizations to obtain the complete list of cyclotomic factors of $F$ with the following lemma:

\begin{lemma}
\label{factome}
Let $m,e\in\mathbb{N}^*$. We set $e_1:=\displaystyle\prod_{p|m}p^{\ord_p(e)}$ and $e_2:=e/e_1$. Then
$$
\Phi_m(X^e)=\Phi_{me_1}(X^{e_2})=\prod_{d|e_2}\Phi_{me_1d}(X).
$$
\end{lemma}
\begin{demo}

To prove the first equality, it suffices to remark that the roots of $\Phi_{me_1}(X^{e_2})$ are also roots of $\Phi_m(X^e)$. Since those polynomials are squarefree, have the same degree and leading coefficient $1$, we conclude that they coincide. A similar argument shows that $\Phi_m(X^e)$ divides the polynomial $\prod_{d|e_2}\Phi_{me_1d}(X)$. Moreover, these polynomials have leading coefficient $1$ and same degree. Indeeed
$$
\sum_{d|e_2}\varphi(d)=e_2,
$$
where $\varphi$ denotes the Euler totient function. Finally these three polynomials are equal.
\end{demo}

In the rest of this section, we establish the following result:
\begin{theorem}
\label{Cyclotomique}
There exists a deterministic algorithm which has the following propriety: for $F\in\Z[X]$ given by its lacunary en\-co\-ding, the algorithm determines a set:
$$
S_F:=\left\{(m_1,e_1),\ldots,(m_t,e_t)\right\}
$$
of couples of integers such that
$$
V(F)_{\rm tors}=\bigcup_{i=1}^tV\left(\Phi_{m_i}\left(X^{e_i}\right)\right).
$$
Furthermore, this algorithm requires
$$
O\left(N^N(M(d)\log\log(d)+h)\right)
$$
bit operations, where $d$ denotes the degree of $F$, $h$ its height and $N$ its number of non zero terms.
\end{theorem}
\vskip 0.5cm

The algorithm underlying Theorem \ref{Cyclotomique} is Algorithm 1 in Subsection 2.2. This result is essentially  \cite[Thm C]{FGS08} but we give a much simpler proof which allows us to explicit the dependence of the complexity on the number of non zero terms $N$ and to speed up its practical execution. More importantly, this proof extends to the multivariate case, as we will see in the next section. We first give a few preliminary results before proving Theorem \ref{Cyclotomique} and finally we will construct a family of polynomials having many <<separated>> cyclotomic factors in Subsection 2.3.

\subsection{Preliminary results}

We first recall the cost of standard arithmetic in $\Z$:

The product of two integers of bit length bounded by $\log(d)$ can be computed in
$$
M(d)=O(\log (d)\log\log (d)\log\log\log (d))=O_{\varepsilon}((\log(d))^{1+\varepsilon})\hskip 0.5cm \forall\,\varepsilon>0
$$
bit operations with the algorithm of Schönhage-Strassen \cite{SS71}.

The gcd of two integers of bit length bounded by $\log(d)$ can be computed in
$$
O(\log (d)(\log\log (d))^2\log\log\log (d))=O_{\varepsilon}((\log(d))^{1+\varepsilon})\hskip 0.5cm \forall\,\varepsilon>0
$$
bit operations with the algorithm of Knuth-Schönhage \cite{K70}.
We refer to \cite{GG03} for the description and analysis of those algorithms.
\vskip 0.5cm
If $\underline{\zeta}\in\mu_{\infty}^n$, we denote by $\ord(\underline{\zeta})$ the order of  $\underline{\zeta}$ in the group $\mu_{\infty}^n$, which equals the less common multiple of the order of its coordinates. For $m\in\N^*$, we also define:
$$
\Psi(m):=2+\sum_{p|m}(p-2),
$$
where $p$ runs over the prime numbers dividing $m$.
The main ingredient of the proof of Theorem \ref{Cyclotomique} is the following result due to Conway and Jones \cite{CJ76}:
\begin{theorem}(Conway-Jones 1976)
\label{CJ}
Let $\zeta_m$ be a root of unity of order $m$. Let also $a_1,\ldots,a_N,$ $\alpha_1,\ldots,\alpha_N$ be integers and $S:=a_1\zeta_m^{\alpha_1}+\cdots+a_N\zeta_m^{\alpha_N}$. If
\begin{enumerate}
\item \label{1}$S=0$,
\item \label{2} no proper subsum of $S$ vanishes,
\item \label{3} $\gcd(\alpha_2-\alpha_1,\ldots,\alpha_N-\alpha_1,m)=1$,
\end{enumerate}
then $m$ is squarefree and $\Psi(m)\leq N$.
\end{theorem}
We say that a vanishing sum of roots of unity is {\it minimal} if the condition \ref{2} is satisfied. Condition \ref{3} is equivalent to
$$
\ord(\zeta_m^{\alpha_2-\alpha_1},\ldots,\zeta_m^{\alpha_N-\alpha_1})=m.
$$
We also need the following lemma:
\begin{lemma}
\label{psiprod}
Let $m_1,\ldots,m_s\in\N^*$, then
$$
\Psi\left(m_1\cdot\cdot\cdot m_s\right)\leq \sum_{i=1}^s \Psi(m_i)-2(s-1).
$$
\end{lemma}
{\bf Proof:} We have that
\begin{align*}
\Psi(m_1\cdots m_{s})&=2+\displaystyle\sum_{p|m_1\cdots m_{s}}(p-2)\\
&\leq 2+\displaystyle\sum_{i=1}^s\displaystyle\sum_{p|m_i}(p-2)
\leq \sum_{i=1}^{s} \Psi(m_i)-2(s-1) \hskip 2.9cm \Box
\end{align*}
\vskip 0.7cm

We also need some upper bounds for the cardinality of certain sets to compute a bound for the complexity of Algorithm $1$:
\begin{lemma}
\label{E_N}
Let $N\in\N^*$ and $E_N$ be the set of partitions of $\{1,\ldots,N\}$ which only contain subsets with at least two elements. Then:
$$
\#E_N\leq N!
$$
\end{lemma}
\begin{demo}

Let $\sigma$ be an element of $\mathfrak{S}_N$, the group of permutations of the set $\{1,\ldots,N\}$. Let $\sigma=c_1\circ\cdots\circ c_s$ be its decomposition into disjoint cycles, then we associate to $\sigma$ the partition $\{J_1,\ldots,J_s\}$ of the set $\{1,\ldots,N\}$ such that $J_j$ is the support of $c_j$ for all  $1\leq j\leq s$. Since this application is a surjection between $\mathfrak{S}_N$ and the set of partitions of the set $\{1,\ldots,N\}$ and since $E$ is strictly included in this set, we have $\#E_N\leq N!$.
\end{demo}
{\bf Remark:} We can improve this estimate to $\#E_N\leq\frac{N!}{e}$ by consi\-dering the subset of $\mathfrak{S}_N$ consisting of the permutations without fixed points, whose image also contains $E_N$.

\begin{lemma}
\label{Q_N}
Let $N,n\in\N$. We set
$$
Q_N:=\{m\in\N \,:\,\Psi(m)\leq N \text{ and m is squarefree}\}
$$
and
$$
Q_{n,N}:=\left\{\underline{\omega}\in\mu_{\infty}^n : \ord(\underline{\omega})\in Q_N\right\}.
$$
Then 
$$
\#Q_N\leq \exp\left(3\sqrt{N\log N}\right)
$$
and
$$
\#Q_{n,N}\leq \exp\left(3(n+1)\sqrt{N\log N}\right).
$$
\end{lemma}
\begin{demo}

We first prove the first inequality. Let $m\in Q_N$ and let $m=\prod_{i=1}^rp_i$ its factorization into primes. Since $\Psi(m)\leq N$, we have
$$
\displaystyle\sum_{i=1}^{r}p_i\leq N+2(r-1).
$$
Furthermore, 
$$
N-2\geq \displaystyle\sum_{i=1}^r(p_i-2)\geq \displaystyle\sum_{i=2}^r(p_i-2)\geq \displaystyle\sum_{i=1}^{r-1}(2i-1)= (r-1)^2.
$$
So, $r\leq \sqrt{N-2}+1$ and
\begin{equation}
\label{r1}
\displaystyle\sum_{i=1}^{r}p_i\leq N+2\sqrt{N-2}.
\end{equation}
Let us now obtain a better upper bound for $r$. By direct computation, we obtain that 
$$
r\leq 3\sqrt{\frac{N}{\log N}}\hskip 1cm {\rm for }\,\, N\leq 10\,000.
$$
We assume now that $N\geq 10\,000$. An effective version of Tchebychev's Theorem (see \cite{RS62}) gives the inequalities:
\begin{equation}
\label{Tchebychev}
\pi(x)\leq 1.26\,\frac{x}{\log x}\hskip 0.5cm \text{and}\hskip 0.5cm \sum_{p\leq x}p\geq \frac{x^2}{2\log x}\hskip 1cm \text{for }x\geq 347.
\end{equation}
Then for $x:=1.15\,\sqrt{N\log N}$, we have $x\geq347$ and thus:
\begin{equation}
\label{r2}
\sum_{p\leq x}p\geq \frac{1.15^2N\log N}{2\log(1.15\,\sqrt{N\log N})}\geq N+2\sqrt{N-2}.
\end{equation}
The last inequality comes from a simple study of function. Combining the two inequalities $(\ref{r1})$ and $(\ref{r2})$, we obtain $\displaystyle\sum_{p\leq x}p\geq\displaystyle\sum_{i=1}^{r}p_i$ and then:
$$
r\leq\pi(1.15\,\sqrt{N\log N})\leq 1.26\,\frac{1.15\,\sqrt{N\log N}}{\log(1.15\,\sqrt{N\log N})}\hskip 0.5cm \text{\rm according to (\ref{Tchebychev}}),
$$
and finally
$$
r\leq 1.26\,\frac{1.15\,\sqrt{N\log N}}{\log(\sqrt{N})}\leq 3\sqrt{N/\log N}.
$$
So each integer in the set $Q_N$ has at most $3\sqrt{N/\log N}$ prime factors which are all smaller than $N$. It follows that
$$
\#Q_N\leq N^{3\,\sqrt{N/\log N}}= \exp\left({3\sqrt{N\log N}}\right) .
$$

Now we prove the second inequality. Let $\underline{\omega}\in Q_{n,N}$ of order $m\in Q_N$. We have established that $m\leq \exp\left(3\sqrt{N\log N}\right)$. Moreover, $\underline{\omega}$ can be represented by the $(n+1)-$uplet~:
$$
(d_1,\ldots,d_n,m),
$$
where the integers $d_i$ are between $0$ and $m-1$ so that:
$$
\underline{\omega}=\left(\zeta_m^{d_1},\ldots,\zeta_m^{d_n}\right),
$$
where $\zeta_m=\exp\left({\frac{2i\pi}{m}}\right)$. Since there is at most $\exp\left(3(n+1)\sqrt{N\log N}\right)$ such $(n+1)-$uplets, the result follows.
\end{demo}

\subsection{Proof of Theorem \ref{Cyclotomique}}

In the following, let 
$$
F(X)=\displaystyle\sum_{i=1}^Na_iX^{\alpha_i}\in\Z[X]
$$
and $\zeta_m$ be a root of unity of order $m$ such that \begin{equation}
\label{cond1}
F(\zeta_m)=0.
\end{equation}
We will determine equivalent conditions to the equation (\ref{cond1}) that we could test algorithmically. Let $\{J_1,\ldots,J_s\}$ be a partition of the set $\{1,\ldots,N\}$ such that 
$$
\displaystyle\sum_{i\in J_j} a_i \zeta_m^{\alpha_i}=0 \hskip 1cm \forall\,\,1\leq j\leq s
$$
and such that every such sum is minimal. For $1\leq j\leq s$, let 
$$
F_j(X):=\displaystyle\sum_{i\in J_j}a_iX^{\alpha_i}
$$
such that $F(X)=\sum_{j=1}^{s}F_j(X)$. We also define for $1\leq j\leq s$, integers $e_j,b_j$ and polynomials $G_j$ such that 
\begin{equation}
\label{Fj}
F_j(X)=X^{b_j}G_j(X^{e_j}),
\end{equation}
where $G_j(0)\neq0$ and such that the exponents of the monomials appearing in $G_j$ are coprime. We finally set:
$$
m_j:=\frac{m}{\gcd(m,e_j)}.
$$
If $\zeta_{m_j}$ denotes a root of unity of order $m_j$, the equation (\ref{cond1}) is then equivalent to:
\begin{equation}
\label{cond2}
G_j(\zeta_{m_j})=0\text{ and } m_j=\frac{m}{\gcd(m,e_j)} \hskip 1 cm \forall \,\,1\leq j\leq s.
\end{equation}
Furthermore, by construction, the sums $G_j(\zeta_{m_j})$ satisfy the conditions of Theorem \ref{CJ} so we have 
$$
\Psi(m_j)\leq N_j \hskip 1 cm \forall \,\,1\leq j\leq s,
$$
where $N_j:=\#J_j$ and $m_j$ is squarefree.\\
We set $e:=\gcd(e_1,\ldots,e_s)$ and, for $1\leq j\leq s$, we also set $e_j':=e_j/e$. Finally we set
$$
m':=\frac{m}{\gcd(m,e)}
$$
and we observe that
$$
\frac{m}{\gcd(m,e_j)}=\frac{m'}{\gcd(m',e_j')} \hskip 1 cm \forall \,\,1\leq j\leq s.
$$
The condition (\ref{cond2}) is then equivalent to:
\begin{equation}
\label{cond3}
G_j(\zeta_{m_j})=0\text{ and } m_j=\frac{m'}{\gcd(m',e_j')} \hskip 1 cm \forall \,\,1\leq j\leq s.
\end{equation}
Moreover, since $m'|{\rm lcm}(m_1,\ldots,m_s)$, we have:
\begin{align*}
\Psi(m')&\leq \Psi({\rm lcm}(m_1,\ldots,m_s))\\
&\leq \Psi(m_1\cdot\cdot\cdot m_s))\\
&\leq \displaystyle\sum_{j=1}^s\Psi(m_j)-2(s-1) \hskip 0.5cm\text{ according to Lemma \ref{psiprod}}\\
&\leq \displaystyle\sum_{j=1}^s N_j-2(s-1)\\
&\leq N-2(s-1)
\end{align*}
and $m'$ is squarefree since it divides the squarefree number ${\rm lcm}(m_1,\ldots,m_s)$. We also remark that after a partition of $\{1,\ldots,N\}$ and an integer $m'$ are fixed, the condition (\ref{cond3}) can be tested since the polynimials $G_j$ and the integers $m_j, e'_j$ only depend on the partition and on $m'$. The algorithm will test the condition (\ref{cond3}) for all possible choices of partitions and for all possible choices of integers $m'$. We can do this as follows:
\vskip 0.5cm
\noindent{\large{\bf Algorithm 1}}

\vskip 0.3cm

\noindent{\bf Input:} A polynomial $F(X)=\sum_{i=1}^Na_iX^{\alpha_i}$ given by its lacunary encoding.\\
{\bf Output:} The set $S_F$ of Theorem \ref{Cyclotomique}.
\begin{enumerate}
\item $S_F\leftarrow \emptyset$.
\item \label{part}For each partition $\{J_1,\ldots,J_s\}$ of $\{1,\ldots,N\}$ such that $\#J_j\geq2$, for all $1\leq j\leq s$ do:
\begin{enumerate}
\item \label{G_j}For all $1\leq j\leq s$, compute the polynomials $G_j$ and the integers $e_j$ associated to the set $J_j$, as in (\ref{Fj}).
\item \label{e_j'}Compute $e:=\gcd(e_1,\ldots,e_s)$ and for $1\leq j\leq s$, compute $e_j':=e_j/e$.
\item \label{M}For each squarefree integer $m'$ satisfying  $\Psi(m')\leq N-2(s-1)$ do:
\begin{enumerate}
\item \label{m_j}For $1\leq j\leq s$, compute $m_j:=\dfrac{m'}{\gcd(m',e_j')}$.
\item \label{test}If $\Phi_{m_j}(X)|G_j(X)$ for all $1\leq j\leq s$, then do\\
$S_F\leftarrow S_F\cup \{(m',e)\}$.
\end{enumerate}
\end{enumerate}
\item Output $S_F$ and terminate.
\end{enumerate}

According to the discussion above, this algorithm correctly outputs the set $S_F$ of Theorem \ref{Cyclotomique}. We will only explain how we can compute the set $Q_N$ of integers considered in step (2.c) and how we can check the divisibility relations in step (2.c.ii).
Let us recall that the set $Q_N$ was defined in Lemma~\ref{Q_N} by
$$
Q_N=\{m\in\N\,|\,m\,\, \text{\rm is squarefree and }\Psi(m)\leq N\}.
$$
We begin by constructing the set $P_N=\{p_1,\ldots,p_t\}$ of prime numbers less than $N$, which can be done with the Eratosthenes sieve in $O(N^2\log N)$ bit operations. Then we compute the subsets $\{q_1,\ldots,q_r\}$ of $P_N$ for which
\begin{equation}
\label{psiq}
2+\sum_{i=1}^r(q_i-2)\leq N.
\end{equation}
Checking the inequality (\ref{psiq}) requires to sum up at most $N$ integers of bit length less than $\log N$ which can be done in $O(N\log N)$ bit operations for a given subset $\{q_1,\ldots,q_r\}$ of $P_N$. Furthermore, $\#P_N<N$ so this set has less than $2^N$ subsets and checking the inequality (\ref{psiq}) for all of them requires $O(2^NN\log N)$ bit operations. Finally, we start with $Q_N$ equal to the empty set and for each subset $\{q_1,\ldots,q_r\}$ of $P_N$ satisfying the inequality $(\ref{psiq})$, we compute the product $\prod_{i=1}^rq_i$ and we add it in $Q_N$. Computing this product requires to perform at most $N$ products of integers of bit length less than $\log N$ and since this product never exceeds $\exp(3\sqrt{N\log N})$ according to Lemma~\ref{Q_N}, this can be done in $O(N^2\log N)$ bit operations. Performing this operation for every set $\{q_1,\ldots,q_s\}$ can be done in $O(2^NN^2\log N)$ bit operations and this is also the numbers of bit operations needed to compute $Q_N$.

\hskip 0.5cm

Let us now explain how we check the divisibility relation $\Phi_{m_j}(X)|G_j(X)$ at step (2.c.ii). For this we use a simple algorithm described by Filaseta and Schinzel in \cite[Theorem 3]{FS04}. It requires the factorization of $m_j$ and is based on the equivalence
$$
\Phi_{m}(X)|G(X)\Longleftrightarrow X^m-1\mid G(X)\times\prod_{p|m}\left(X^{m/p}-1\right),
$$
where the product runs over the prime numbers dividing $m$. We refer to \cite{FS04} for the complete description of this algorithm and its running time. In our case, since $m_j$ is squarefree, the first step of their algorithm can be avoided and the divisibility relation $\Phi_{m_j}(X)|G_j(X)$ can be checked in $O\left(N^3(\log N)^22^{3\sqrt{N\log N}}(h+M(d))\right)$ bit operations.

\vskip 0.5cm

Finally we estimate the overall complexity of Algorithm $1$. We recall that at step (2.a), the polynomials $G_j$ are defined by:
$$
F_j(X)=X^{b_j}G_j(X^{e_j}).
$$
In order to compute these polynomials, we have to perform $N$ substractions of integers of bit length bounded by $\log(d)$, then we have to perform at most $N$ gcd computations of integers of bit length bounded by $\log(d)$, and finally compute $N$ divisions of integers of bit length bounded by $\log(d)$ and all of these computations require $O(NM(d)\log\log(d))$ bit operations. 

At step (2b), we have to compute $s-1$ gcd of integers of bit length bounded by $\log(d)$ to compute $e$, then  we have to perform $s$ Euclidean divisions of integers of bit length bounded by $\log(d)$ to compute the integers $e_j'$. So in the whole, this requires $O(NM(d)\log\log(d))$ bit operations. 

At step (2c), in order to compute the integers $m_j$, we perform $s$ divisions and $s$ gcd computations of integers of bit length bounded by $d$ and this requires again $O(NM(d)\log\log(d))$ bit operations. In order to check the divisibility relations, we apply $s$ times the algorithm of Filaseta and Schinzel \cite[Theorem 3]{FS04} and this can be done in 
$$
O\left(N^4(\log N)^22^{3\sqrt{N\log N}}(h+M(d))\right)
$$
bit operations.

Moreover, since there are at most $N!$ partitions to consider at step (2) by Lemma \ref{E_N} and since there are at most $\exp\left(3\sqrt{N\log N}\right)$ integers to consider at step (2.c) by Lemma \ref{Q_N}, the total number of bit operations performed by the algorithm is:
$$
O\left(N^N(M(d)\log\log(d)+h)\right),
$$
this terminates the proof of Theorem \ref{Cyclotomique}.\hfill\mbox{$\Box$}

\subsection{A family of polynomials having many <<separated>> cyclotomic factors}
\label{exemple}

We believe that the output of Algorithm $1$ can be exponential in the size of the input in the worst case. The following example, which has been elaborated with Denis Simon, support this conjecture. When $N=2n$ is an even integer, we construct polynomials for which the output of the algorithm can not be regrouped in less than $n!$ pairs of integers.
\begin{proposition}
\label{propexample}
Let $n\in\N$ and $p_1,\ldots,p_{n!}$ be $n!$ prime numbers greater than $2n$. Then there exists a polynomial $f\in\Z[X]$ having $2n$ non zero terms and satisfying
\begin{enumerate}
\item $X^{p_r}-1|f(X)$\hskip 1.85cm $\forall\,\,1\leq r\leq n!$
\item $X^{p_rp_{r'}}-1\nmid f(X)$ \hskip 1cm $\forall\,\, 1\leq r<r'\leq n!$
\end{enumerate}
\end{proposition}
\begin{demo}

\noindent Let $f(X):=\sum_{i=1}^nX^{\alpha_i}-\sum_{i=n+1}^{2n}X^{\alpha_i}$. We will compute values of $\alpha_i$ for which $f$ satisfies the proposition \ref{propexample}.
We consider all the partitions $\{J_{1},\ldots,J_{n}\}$ of the set $\{1,\ldots,2n\}$ made by 2-elements sets $J_1,\ldots,J_n$ such that each $J_j$ contains an element in $\{1,\ldots,n\}$ and an other one in $\{n+1,\ldots,2n\}$. For instance,
$$
\left\{\{1,n+1\},\{2,n+2\},\ldots,\{n,2n\}\right\}
$$
is one of these partitions. We remark that there is $n!$ such partitions. We enumerate these partitions as
$$
\pi_r:=\left\{J_{r,1},\ldots,J_{r,n}\right\}\hskip 1cm\forall\,\,1\leq r\leq n!.
$$
To each partition $\pi_r=\{J_{r,1},\ldots,J_{r,n}\}$ we associate injectively a prime number $p_r$ strictly larger than $2n$. According to the Chinese reminder theorem, we can choose $\alpha_i$ such that:
$$
i\in J_{r,k} \Longleftrightarrow \alpha_i\equiv k-1 \mod p_r \hskip 1cm \forall\,\,1\leq r\leq n!.
$$
By construction, all the $\alpha_i$ are distinct and we have:
$$
X^{p_r}-1\big\vert f(X)\hskip 1cm \forall\,\,1\leq r\leq n!.
$$
We finally prove the second part of the proposition. Let us suppose, for example,  that $X^{p_1p_2}-1$ divides $f(X)$. We then have $\Phi_{p_1p_2}(X)\vert f(X)$. Thus $f(\zeta)=0$, where $\zeta$ denotes a root of unity of order $p_1p_2$. Let also
$$
f(\zeta)=\sum_{i=1}^n\zeta^{\alpha_i}-\sum_{i=n+1}^{2n}\zeta^{\alpha_i}=:S_1+\cdots+S_t,
$$
be a decomposition of $f(\zeta)$ into minimal vanishing sum. Let $m_j$ be the order of $S_j$. Since $m_j|p_1p_2$, we have $m_j\in\{1,p_1,p_2,p_1p_2\}$ and $m_j$ satisfies $\Psi(m_j)\leq \#S_j\leq 2n$ according to Theorem \ref{CJ}. However $m_j=1$ since $\Psi(p_1)=p_1>2n$, $\Psi(p_2)=p_2>2n$ and $\Psi(p_1p_2)=p_1+p_2-2>2n$.So each $S_j$ is a multiple of a vanishing sum of roots of unity of the shape $\pm1\cdots\pm1=0$. Since $S_j$ is minimal, it must have length $2$ and we also have $t=n$. Thus for $1\leq j\leq n$, we can rewrite $S_j$ as:
$$
S_j=\zeta^{\alpha_{a_j}}-\zeta^{\alpha_{b_j}}\
$$
with
$$
1\leq a_j\leq n<b_j\leq 2n\hskip 0.5cm \text{and}\hskip 0.5cm  \alpha_{a_j}\equiv \alpha_{b_j} \mod p_1p_2.
$$
Then
$$
\alpha_{j,1}\equiv\alpha_{j,2}\mod p_1\hskip 0.5cm \text{and}\hskip 0.5cm \alpha_{j,1}\equiv\alpha_{j,2}\mod p_2\hskip 0.5cm \forall\,\,1\leq j\leq n.
$$
Thus we necessarily have $\pi_r=\pi_1$ and $\pi_r=\pi_2$. This is a contradiction since $\pi_1\neq\pi_2$. Thus, $X^{p_1p_2}-1\nmid f(X)$.
\end{demo}

Let $N=2n$ be an even integer and $f$ be the polynomial constructed in the proof of Proposition \ref{propexample}. It is easy to see that the algorithm $1$, when applied to $f$, will output at least $n!$ pairs of integers. Unfortunately, this example do not suffice to prove that the size of the output is not polynomial in the input size, because the degree of $f$ is doubly exponential in $N$. So the size of the input is also exponential in $N$. On the other hand, it proves that for a given number $N$ of non zero terms, the output can be exponential in~$N$.
\begin{example}
{\rm for} $N=6$.
\end{example}
Let
$$
f(x)=x^{\alpha_1}+x^{\alpha_2}+x^{\alpha_3}-x^{\alpha_4}-x^{\alpha_5}-x^{\alpha_6}.
$$
We make the following choices of prime numbers for the $n!=6$ partitions:
\begin{align*}
\{\{1,4\},\{2,5\},\{3,6\}\}\longmapsto & 7\\
\{\{1,4\},\{2,6\},\{3,5\}\}\longmapsto & 11\\
\{\{1,5\},\{2,4\},\{3,6\}\}\longmapsto & 13\\
\{\{1,5\},\{2,6\},\{3,4\}\}\longmapsto & 17\\
\{\{1,6\},\{2,4\},\{3,5\}\}\longmapsto & 19\\
\{\{1,6\},\{2,5\},\{3,4\}\}\longmapsto & 23
\end{align*}
and we consider the following system of congruences:
\begin{center}

\begin{tabular}{|c|c|c|c|c|c|c|}
\hline  &$\mod 7$  &$\mod 11$  &$\mod 13$  &$\mod 17$  &$\mod 19$ &$\mod23$ \\ 
\hline  $\alpha_1$& $0$ & $0$ & $0$ & $0$ & $0$ & $0$ \\ 
\hline  $\alpha_2$& $1$ & $1$ & $1$ & $1$ & $1$ & $1$ \\ 
\hline  $\alpha_3$& $2$ & $2$ & $2$ & $2$ & $2$ & $2$ \\ 
\hline  $\alpha_4$& $0$ & $0$ & $1$ & $2$ & $1$ & $2$ \\ 
\hline  $\alpha_5$& $1$ & $2$ & $0$ & $0$ & $2$ & $1$ \\ 
\hline  $\alpha_6$& $2$ & $1$ & $2$ & $1$ & $0$ & $0$ \\ 
\hline 
\end{tabular} 
\end{center}
A solution of this system is:
\begin{align*}
\alpha_1=&\,\, 0\\
\alpha_2=& \,\,1\\
\alpha_3=& \,\,2\\
\alpha_4=& \,\,6\,282\,199\\
\alpha_5=& \,\,2\,501\,941\\
\alpha_6=& \,\,6\,088\,721
\end{align*}
which gives the polynomial:
$$
f(x)=1+X+X^2-X^{2501941}-X^{6088721}-X^{6282199}
$$
for which the algorithm $1$ outputs the set of pairs:
$$
\{(2,2),\,(1,7),\,(1,11),\,(1,13),\,(1,17),\,(1,19),\,(1,23)\}.
$$
and we can not regroup the pairs, according to Proposition \ref{propexample}.
\begin{example} {\rm for} $N=8$.
\end{example}
We can construct the following polynomial
\begin{align*}
f(x):=&1+x+x^2+x^3-x^{10649315971896428139150081202897150286932}\\
&-x^{10417696267221214855704118228748809801239}\\
&-x^{2747750133111287905524860455880456232062}\\
&-x^{626914938199634951807585972855218426387}
\end{align*}
by considering the $24$ prime numbers between $11$ and $107$.
For this polynomial, Algorithm $1$ outputs the set of pairs:
\begin{align*}
\{&(1,13),\,(1,17),\,(1,19),\,(1,22),\,(1,23),\,(1,31),\,(1,37),\,(1,41),\,(1,43),\\
&\,(1,47),\,(1,53),\,(1,58),\,(1,59),\,(1,61),\,(1,71),\,(1,79),\,(1,83),\,(1,89),\\
&(1,97),\,(1,101),\,(1,103),\,(1,107),\,(1,134),\,(1,146)\}
\end{align*}
which corresponds to $29$ distinct cyclotomic factors and by Proposition \ref{propexample}, these can not be regrouped in less than $24$ pairs, as above.

\section{Torsion cosets of a hypersurface}

\label{Sec3}
We want to compute a representation of the torsion points lying in a hypersurface defined by a polynomial with integer coefficients. In this section, we prove Theorem \ref{variete} in the case where $V$ is a hypersurface and describe the algorithm underlying this theorem in Subsection \ref{algo2}.

\subsection{Preliminary results}
First, we will need the following extension of Theorem \ref{CJ} to the multivariate case.
\begin{corollary}
\label{CJ2}
Let  $\underline{\zeta}\in\mu_{\infty}^n$ be a point of order $m$ and $F(\underline{X})=\sum_{i=1}^Na_i\underline{X}^{\underline{\alpha_i}}\in\Z[X_1^{\pm 1},\ldots,X_n^{\pm 1}]$ be a Laurent polynomial. If
\begin{enumerate}
\item $F(\underline{\zeta})=0$, 
\item no proper subsum of $F(\underline{\zeta})$ vanishes,
\item $\displaystyle\sum_{1\leq i,j\leq N}^N\Z(\underline{\alpha_i}-\underline{\alpha_j})=\Z^n$,
\end{enumerate}
then $m$ is squarefree and $\Psi(m)\leq N$.

\end{corollary}
\begin{demo}

We just have to prove that $S:=\sum_{i=1}^Na_i\underline{\zeta}^{\underline{\alpha_i}-\underline{\alpha_1}}$ satisfies the three conditions of Theorem \ref{CJ}. The first two conditions are satisfied since $F(\underline{\zeta})$ is a minimal vanishing sum by hypothesis. Let $M:=\underset{1\leq i\leq N}{{\rm lcm}}\left(\ord(\underline{\zeta}^{\underline{\alpha_i}-\underline{\alpha_1}})\right)$. It remains to check the last condition, i.e. to check if  $M=m.$ Since the divisibility $M|m$ is clear, we only have to prove the reverse divisibility condition. In the following, we set $\underline{\zeta}=(\zeta_m^{d_1},\ldots,\zeta_m^{d_n})$, where $\zeta_m$ is a root of unity of order $m$. Let us first remark that for all  $b_1,\ldots,b_N\in\Z$, we have:
$$
\ord\left(\underline{\zeta}^{\sum_{i=1}^Nb_i(\underline{\alpha_i}-\underline{\alpha_1})}\right)\bigg\vert M.
$$
Moreover, since $\displaystyle\sum_{i=1}^N\Z(\underline{\alpha_i}-\underline{\alpha_1})=\Z^n$, there exists $b_1,\ldots,b_N\in\Z$ such that\\ $\displaystyle\sum_{i=1}^Nb_i(\underline{\alpha_i}-\underline{\alpha_1})=(1,0,\ldots,0)$. Thus,
$$
\ord\left(\zeta_m^{d_1}\right)=\ord\left(\underline{\zeta}^{\sum_{i=1}^Nb_i(\underline{\alpha_i}-\underline{\alpha_1})}\right)\bigg\vert M
$$
and similarly
$$
\ord\left(\zeta_m^{d_i}\right)\big\vert M\hskip 1cm \text{pour }2\leq i\leq N.
$$
Finally we have 
$$
m=\underset{1\leq i\leq N}{{\rm lcm}}\left(\ord\left(\zeta_m^{d_i}\right) \right)\big\vert M.
$$
and the last condition of Theorem \ref{CJ} is satisfied.
\end{demo}

We will also need to be able to compute a basis of a subgroup of $\Z^n$ from a generating family of vectors. Let $\underline{\beta_1},\ldots,\underline{\beta_r}$ $\in\Z^n$ be vectors generating a subgroup $R$ of $\Z^n$ of rank $k\leq r$ and $M$ the matrix of size $r\times n$ whose rows are the vectors $\underline{\beta_1},\ldots,\underline{\beta_r}$. We can compute from $M$ a matrix $H$ of size $k\times n$, which is upper triangular, whose rows consist of a basis of $R$. The matrix $H$ is called the {\it Hermite normal form} of $M$ and can be computed with the following lemma:
\begin{lemma}
\label{HNF}
Let $n,d,k$ be integers, $M$ a matrix of size $k\times n$ with integer coefficients bounded by $d$ in absolute value. Then we can compute $H$, the Hermite normal form of $M$, in $O\left(k^4n^2M(d)\log\log(d)\right)$ bit operations.
\end{lemma}
In \cite{LS96}, the authors give a proof of this result with a better complexity but this one will suffice for our need. We refer to \cite{GG03} for more details about the Hermite normal form of an integral matrix.

\subsection{Proof of Theorem \ref{variete} for a hypersurface}
\label{algo2}
Let $F(\underline{X})=\sum_{i=1}^Na_i\underline{X}^{\underline{\alpha_i}}\in\Z[X_1,\ldots,X_n]$, and $\underline{\zeta}\in\mu_{\infty}^n$  such that
\begin{equation}
\label{Cond1}
F(\underline{\zeta})=0.
\end{equation}
As we have done in the proof of Theorem \ref{Cyclotomique}, we will find some conditions equivalent to (\ref{Cond1}) that we could test algorithmically.

Let $\{J_1,\ldots,J_s\}$ be a partition of $\{1,\ldots,N\}$ such that: 
$$
\sum_{i\in J_j} a_i \underline{\zeta}^{\underline{\alpha_i}}=0\hskip 1cm \forall\,\,1\leq j\leq s
$$
and such that every of these sums is minimal. We set $N_j:=\#J_j$ and, for $1\leq j\leq s$,
$$
F_j(\underline{X}):=\displaystyle\sum_{i\in J_j}a_i\underline{X}^{\underline{\alpha_i}}.
$$
Thus we have $F(\underline{X})=\sum_{j=1}^{s}F_j(\underline{X})$. Then we renumber the coefficients and the exponents of $F_j$ in the following manner:
$$
F_j(\underline{X}):=\displaystyle\sum_{i=1}^{N_j}a_{j,i}\underline{X}^{\underline{\alpha_{j,i}}}\,\,.
$$
Let, for $1\leq j\leq s$:
$$
R_j:=\displaystyle\sum_{i=1}^{N_j}\Z(\underline{\alpha_{j,1}}-\underline{\alpha_{j,i}})
$$
be the subgroup of $\Z^n$ spanned by the differences of exponents of $F_j$. We denote by $k_j$ the rank of this group and by $\underline{\lambda_{j,1}},\ldots,\underline{\lambda_{j,k_j}}$ a basis of $R_j$. Then we can find a Laurent polynomial $G_j$ in $k_j$ variables such that:
\begin{equation}
\label{F_j2}
F_j(\underline{X})=\underline{X}^{\underline{\alpha_{j,1}}}G_j(\underline{X}^{\underline{\lambda_{j,1}}},\ldots,\underline{X}^{\underline{\lambda_{j,k_j}}}).
\end{equation}
We set for $1\leq l\leq k_j$, 
$$
\omega_{j,l}:=\underline{\zeta}^{\underline{\lambda_{j,l}}}.
$$
We have
\begin{equation}
\label{G_j=0}
G_j(\omega_{j,1},\ldots,\omega_{j,k_j})=0, \hskip 1cm {\rm for }\,\,1\leq j \leq s.
\end{equation}
The Laurent polynomial $G_j$ and the point $(\omega_{j,1},\ldots,\omega_{j,k_j})$ of $\mu_{\infty}^{k_j}$ satisfy the conditions of Corollary \ref{CJ2}. If $m_j$ denotes the order of $(\omega_{j,1},\ldots,\omega_{j,k_j})$ in $\mu_{\infty}^{k_j}$, we conclude that $m_j$ is squarefree and $\Psi(m_j)\leq N_j$. Now let $R:=\sum_{j=1}^sR_j$  and $\underline{\lambda_1},\ldots,\underline{\lambda_k}$ be a basis of $R$. For $1\leq t\leq k$,  we set $\omega_t:=\underline{\zeta}^{\underline{\lambda_t}}$ and $m':=\ord(\omega_1,\ldots,\omega_k)$. By construction, we have:
$$
m'|{\rm lcm}(m_1,\ldots,m_s).
$$
As in the proof of Theorem \ref{Cyclotomique}, the integer $m'$ is squarefree and satisfies $\Psi(m')\leq N-2(s-1)$.\\
If the equalities (\ref{G_j=0}) are satisfied for $1\leq j\leq s$, we have $F(\underline{\zeta})=0$ for all $\underline{\zeta}$ satisfying:
$$
\omega_t=\underline{\zeta}^{\underline{\lambda_t}}\hskip 1cm\text{for }1\leq t\leq k.
$$
Thus we have
$$
F(\underline{\zeta})=0\hskip 1cm \forall\,\,\underline{\zeta}\in B\left(L,(\omega_1,\ldots,\omega_k)\right),
$$
where $L$ denotes the matrix of size $k\times n$ whose rows are the vectors $\underline{\lambda_1},\ldots,\underline{\lambda_k}$.
We recall that $B(L,(\omega_1,\ldots,\omega_k))$ was defined in (\ref{eqoutput}) as
$$
B(L,(\omega_1,\ldots,\omega_k))=\left\{\underline{x}\in\G_m^n\mid\underline{x}^L=(\omega_1,\ldots,\omega_k)\right\}.
$$

Let us summarize what we just proved. If $F(\underline{\zeta})=0$ for some $\underline{\zeta}\in\mu_{\infty}^n$, then there exist $L\in M_{n,k}(\Z)$ and $\underline{\omega}\in\mu_{\infty}^k$, where $k\leq n$, such that
$$
\underline{\zeta}\in B(L,\underline{\omega})\hskip 0.5cm {\rm and}\hskip 0.5cm B(L,\underline{\omega})\subset V(F).
$$
Furthermore, the matrix $L$ corresponds to a partition of the support of $F$ and can be computed from it. On the other hand
$$
\underline{\omega}\in Q_{k,N}=\left\{\underline{\zeta}\in\mu_{\infty}^k\,:\,\Psi(\ord(\underline{\zeta}))\leq N\right\},
$$
which has cardinality bounded by $\exp(3(k+1)\sqrt{N\log N})$ by Lemma \ref{Q_N}.
Thus the algorithm will test the equa\-lities (\ref{G_j=0}) for all possible partitions of the support of $F$ and for all possible choices of points in $Q_{n,N}$. We can do it as follows:

\vskip 0.5cm

\noindent{\large\bf Algorithm 2}

\vskip 0.3cm

\noindent{\bf Input:} A polynomial $F(\underline{X})=\sum_{i=1}^Na_i\underline{X}^{\underline{\alpha_i}}$ given by its lacunary encoding.\\
{\bf Output:} A list $S_F$ of torsion cosets satisfying Theorem \ref{variete}.
\begin{enumerate}
\item \label{E1}$S_F\leftarrow \emptyset$.
\item \label{E2}For each partition $\{J_1,\ldots,J_s\}$ of $\{1,\ldots,N\}$ such that $\#J_j\geq2$ for all $1\leq j\leq s$, do:
\begin{enumerate}

\item \label{E2a1}Compute for $1\leq j\leq s$, the Laurent polynomials $G_j$ and the vectors $\underline{\lambda_{j,1}},\ldots,\underline{\lambda_{j,k_j}}$ as in (\ref{F_j2}).

\item \label{E2a2}Compute a basis $\underline{\lambda_1},\ldots,\underline{\lambda_k}$ of $R$ and compute, for $1\leq j\leq s$ and for $1\leq h\leq k_j$ the integers $\delta_{j,h,t}$ such that $\underline{\lambda_{j,h}}=\sum_{t=1}^k\delta_{j,h,t}\underline{\lambda_t}$.

\item \label{E2a}For each squarefree integer $m'$ such that $\Psi(m')\leq N-2(s-1)$ and for each $(\omega_1,\ldots,\omega_k)\in\mu_{\infty}^k$ of order $m'$ do:
\begin{enumerate}

\item \label{E2a3}Compute for $1\leq j\leq s$ and for $1\leq h\leq k_j$, $\omega_{j,h}:=\prod_{t=1}^k\omega_t^{\delta_{j,h,t}}$ and set $m_j:=\ord(\omega_{j,1},\ldots,\omega_{j,k_j})$.

\item \label{E2a4} If for all $1\leq j\leq s$, we have $G_j(\omega_{j,1},\ldots,\omega_{j,k_j})=0$, then do $S_F\leftarrow S_F\cup \left\{\left(L,(\omega_1,\ldots,\omega_k\right))\right\}$, where $L$ is the matrix of size $k\times n$ whose rows are the vectors $\underline{\lambda_1},\ldots,\underline{\lambda_k}$.

\end{enumerate}

\end{enumerate}
\item Output $S_F$ and terminate.
\end{enumerate}

The discussion above the description of Algorithm $2$ guarantee its correctness and it only remains to evaluate its running time.

At step (2.a) we first have to compute, for $1\leq j\leq s$, the differences $\underline{\alpha_{j,i}}-\underline{\alpha_{j,1}}$. This requires $O(Nn\log(d))$ bit operations. Then we have to compute a basis $(\underline{\lambda_{j,1}},\ldots,\underline{\lambda_{j,k_j}})$ of the subgroup $R_j$ of $\Z^n$. This can be done by computing the Hermite normal form of the matrix whose rows consist of the vectors  $\underline{\alpha_{j,i}}-\underline{\alpha_{j,1}}$. This requires $O(N^4n^2M(d)\log\log(d))$ bit operations by Lemma \ref{HNF}. In order to compute the Laurent polynomials $G_j$, we then have to compute the coordinates of the vectors  $\underline{\alpha_{j,i}}-\underline{\alpha_{j,1}}$ in the basis $(\underline{\lambda_{j,1}},\ldots,\underline{\lambda_{j,k_j}})$. Thus we have to invert a linear system of size $n\times n$ with coefficients bounded by $2d$. This can be done in $O(n^3M(d)\log \log(d))$ bit operations. Finally, we can bound the total number of bit operations performed at this step by
$$
O\left(N^5n^3M(d)\log\log(d)\right).
$$

At step (2.b), we compute the Hermite normal form of the matrix whose rows are the vectors $\underline{\lambda_{j,h}}$, where $1\leq j\leq s$ and $1\leq h\leq k_j$. We obtain a basis $(\underline{\lambda_1},\ldots,\underline{\lambda_k})$ of $R$. Finally we compute the coordinates of the vectors $\underline{\lambda_{j,h}}$ in this basis. This step again requires
$$
O\left(N^5n^3M(d)\log\log(d)\right)
$$
bit operations.

At step (2.c), each number $\omega_{t}$ is represented by the pair of integers $(\rho_t,m')$ such that
$$
\omega_t=\exp\left(\frac{2i\rho_t\pi}{m'}\right).
$$

At step (2.c.i),  each number $\omega_{j,h}$ is represented by the pair of integers
$$
\left(d_{j,h},m'\right):=\left(\sum_{t=1}^k\delta_{j,h,t}\rho_t \mod m',m'\right)
$$
such that
$$
\omega_{j,h}=\exp\left(\frac{2id_{j,h}\pi}{m'}\right).
$$
This computation requires 
$$
O\left(nM(d)M\left(\exp\left(3\sqrt{N\log N}\right)\right)\right)=O\left(nN^2M(d)\right)
$$
bit operations, since $m'\leq \exp\left(3\sqrt{N\log N}\right)$ according to Lemma \ref{Q_N}. Then we can compute for $1\leq j\leq s$, the integer
$$
m_j:=\underset{1\leq h\leq k_j}{{\rm lcm}}\left(\dfrac{m'}{\gcd(m',d_{j,h})}\right).
$$
These computations require $O(nN^2)$ bit operations. Finally this step requires
$$
O(nN^2M(d))
$$
bit operations.

At step (2.c.ii), we have to test, for $1\leq j\leq s$ the equatlity
$$
G_j(\omega_{j,1},\ldots,\omega_{j,k_j})=0.
$$
This is equivalent to test, for $1\leq j\leq s$ the divisibility condition
$$
\Phi_{m_j}(X)\mid G_j(X^{d_{j,1}},\ldots,X^{d_{j,k_j}})
$$
which can be checked by the method used in the step (2.c.ii) of Algorithm $1$ in
$$
O\left(N^4(\log N)^2 2^{3\sqrt{N\log N}}(h+M(d))\right)
$$
bit operations.

Moreover, since there are at most $N!$ partitions to consider at step (2) according to Lemma \ref{E_N} and since there are at most $\exp\left(3(n+1)\sqrt{N\log N}\right)$ torsion points to consider at step (2.c) by Lemma \ref{Q_N}, the total number of bit operations performed by the algorithm is:
$$
O\left( N^{nN}(M(d)\log\log(d)+h)\right).
$$
This terminates the proof of Theorem \ref{variete} in the case of a hypersurface.

\section{Computing the torsion cosets of a variety}

\label{Sec4}
The aim of this section is to prove Theorem \ref{variete}. In the following, we fix $F_1,\ldots,F_k\in\Z [X_1,\ldots,X_n]$ and we would like to compute a representation of the torsion cosets included in the variety:
$$
V:=\left\{ \underline{x}\in\G_m^n : F_1(\underline{x})=\ldots=F_k(\underline{x} ) \right\}.
$$
We will explain how we can modify Algorithm $2$ of the last section to treat the case of several polynomials. We remark that the same modifications can be applied to Algorithm $1$ to compute the common cyclotomic factors of several univariate polynomials.

Let $\underline{\zeta}\in V\cap\mu_{\infty}^n$ and $N_i$ be the number of non zero terms in $F_i$, for $1\leq i\leq k$ and let also $N:=\max(N_1,\ldots,N_k)$. For each integer $i$ between $1$ and $k$, we can reduce the equality $F_i(\underline{\zeta})=0$ to minimal sums as we have done in previous sections, by considering a partition of the set $\{1,\ldots,N_i\}$. Once performed this step for every integer $i$, we are reduced to step (\ref{E2}) of Algorithm $2$ and all other steps of Algorithm $2$ can be kept.\\
The number of partitions to consider can be bounded now by
$$
N_1!\times\cdots \times N_k!\leq (N!)^k.
$$
The number of torsion points which appear at step (2.c) can be bound now by
$$
\exp\left({3(n+1)\sqrt{kN\log kN}}\right).
$$
Finally, the total number of bit operations executed by the algorithm is:
$$
O\left(N^{nkN}\left(M(d)\log\log(d)+h\right) \right).
$$
This concludes the proof of Theorem \ref{variete}.\hfill\mbox{$\Box$}

\section{An upper bound for the number of torsion cosets}

We prove Corollary \ref{cor}. Let $V$ be a variety defined by $k$ polynomials, each having at most $N$ non zero terms. Our algorithm outputs a representation of $\overline{V_{\rm tors}}$ as the union of torsion cosets. We remark that at most one torsion coset is output in every loop performed by Algorithm $2$. Thus the total number $t$ of torsion cosets included in $V$ is bounded by the number of loops performed by the algorithm. Every loop corresponds to a partition of these polynomials and to a choice of a point in $\mu_{\infty}^n$ of squarefree order $m$ satisfying $\Psi(m)\leq N$. As we have seen before, there are at most $(N!)^k$ partitions and the number of points in $\mu_{\infty}^n$ to consider is bounded by $\exp\left({3(n+1)\sqrt{kN\log kN}}\right)$. Thus
$$
t\leq (N!)^k\exp\left({3(n+1)\sqrt{kN\log kN}}\right),
$$ 
which is the desired upper bound.\hfill\mbox{$\Box$}

\end{document}